\documentclass[reqno, a4paper, 10pt]{amsart}
\usepackage{amssymb,url, color, mathrsfs, multirow, makecell, array, caption}
\usepackage[colorlinks=true, bookmarks=true, pdfstartview=FitH, pagebackref=true, linktocpage=true]{hyperref}
\usepackage[short,nodayofweek]{datetime}
\usepackage[pagewise,running,mathlines,displaymath, switch]{lineno}
\usepackage{times}
\usepackage{indentfirst, tabularx}
\usepackage[table]{xcolor}
\usepackage{float, hhline, colortbl}
\usepackage{graphicx}
\usepackage{longtable}
\usepackage{empheq}
\usepackage{aligned-overset}

\usepackage[framemethod=TikZ]{mdframed}
 
\surroundwithmdframed[
 topline=false,
 rightline=true,
 bottomline=false,
 leftmargin=\parindent,
 skipabove=\medskipamount,
 skipbelow=\medskipamount
]{theorem}

\theoremstyle{plain}
\numberwithin{equation}{section}
\newtheorem{theorem}{Theorem}[section]
\newtheorem{lemma}[theorem]{Lemma}

\theoremstyle{remark}

\newtheorem{remark}[theorem]{Remark}

\DeclareMathOperator{\R}{\mathbf{R}}

\hypersetup{
pdftitle={All conditions for Stein-Weiss inequalities are necessary},
pdfauthor={Qu\^oc Anh Ng\^o},
colorlinks = true,
linkcolor = magenta,
citecolor = blue,
}

\def\R{\mathbf R}

\numberwithin{equation}{section}

\makeatletter
\def\@cite#1#2{[\textbf{#1}\if@tempswa, #2\fi]}
\makeatother

\title[All conditions for Stein--Weiss inequalities are necessary]{All conditions for Stein--Weiss inequalities are necessary}

\def\cfac#1{\ifmmode\setbox7\hbox{$\accent"5E#1$}\else\setbox7\hbox{\accent"5E#1}\penalty 10000\relax\fi\raise 1\ht7\hbox{\lower1.05ex\hbox to 1\wd7{\hss\accent"13\hss}}\penalty 10000\hskip-1\wd7\penalty 10000\box7 }

\author[Q.A. Ng\^o]{Qu\cfac oc Anh Ng\^o}
\address[Q.A. Ng\^o]{
University of Science, Vietnam National University, Hanoi, Vietnam
\\
ORCID iD: 0000-0002-3550-9689}
\email{\href{mailto: Q.A. Ng\^o <nqanh@vnu.edu.vn>}{nqanh@vnu.edu.vn}
}

\allowdisplaybreaks

\begin{document}

\begin{abstract}
The famous Stein--Weiss inequality on $\R^n \times \R^n$, also known as the doubly weighted Hardy--
Littlewood--Sobolev inequality, asserts that
\[
\Big| \iint_{\R^n \times \R^n} \frac{f(x) g(y)}{|x|^\alpha |x-y|^\lambda |y|^\beta} dx dy \Big| 
\lesssim \| f \| _{L^p(\R^n)} \| g\| _{L^r(\R^n)}
\]
holds for any $f\in L^p(\R^n)$ and $g\in L^r(\R^n)$ under several conditions on the parameters $n$, $p$, $r$, $\alpha$, $\beta$, and $\lambda$. Extending the above inequality to either different domains rather than $\R^n \times \R^n$ or classes of more general kernels rather than the classical singular kernel $|x-y|^{-\lambda}$ has been the subject of intensive studies over the last three decades. For example, Stein--Weiss inequalities on the upper half space, on the Heisenberg group, on homogeneous Lie group are known. Served as the first step, this work belongs to a set in which the following inequality on the product $\R^{n-k} \times \R^n$ is studied
\[
\Big| \iint_{\R^n \times \R^{n-k}} \frac{f(x) g(y)}{|x|^\alpha |x-y|^\lambda |y|^\beta} dx dy \Big| 
\lesssim \| f \| _{L^p(\R^{n-k})} \| g\| _{L^r(\R^n)}.
\]
Toward the validity of the above new inequality, in this work, by constructing suitable counter-examples, we establish all conditions for the parameters $n$, $p$, $r$, $\alpha$, $\beta$, and $\lambda$ \textit{necessarily} for the validity of the above proposed inequality. Surprisingly, these necessary conditions applied to the case $k=1$ suggest that the existing Stein--Weiss inequalities on the upper half space are yet in the optimal range of the parameter $\lambda$. This could reflect limitations of the methods often used. Comments on the Stein--Weiss inequality on homogeneous Lie groups as well as the reverse form for Stein--Weiss inequalities are also made.
\end{abstract}

\date{\bf \today \, at \currenttime}



\maketitle


\section{Introduction}

The classical Hardy--Littlewood--Sobolev inequality on $\R^n$ states that for any $n\geq 1$, $p,r > 1$, and $\lambda \in (0,n)$ satisfying the so-called balance condition
\begin{equation*}\label{eq-HLS-BC}
1 /p + 1 /r +\lambda /n =2,
\end{equation*}
the following convolution-type inequality holds
\begin{equation}\label{eq-HLS}
\Big |\iint_{\R^n \times \R^n} \frac{f(x)g(y)}{|x-y|^\lambda} dx dy\Big| 
\lesssim \| f \| _{L^p(\R^n )} \| g\| _{L^r(\R^n)}
\end{equation}
for any $f\in L^p(\R^n)$ and any $g\in L^r(\R^n)$. The case $n=1$ of \eqref{eq-HLS} was first appeared in a work due to Hardy and Littlewood in 1928; see \cite{hl1928}. Then, in 1938 Sobolev generalized the inequality \eqref{eq-HLS} for arbitrary $n$ in \cite{sobolev1938}. The sharp form of \eqref{eq-HLS} as well as the existence of a pair of the optimal functions was eventually obtained by Lieb in \cite{Lieb} via a symmetric decreasing rearrangement technique; see also \cite{Lions} for a different argument based on concentration-compactness principle.

In 1958, Stein and Weiss proved the so-called doubly weighted Hardy--Littlewood--Sobolev inequality, which generalizes the HLS inequality \eqref{eq-HLS} by inserting two weights $|x|^{-\alpha}$ and $|y|^{-\beta}$ into the integrand $f(x)|x-y|^{-\lambda}g(y)$. To be more precise, it was proved in \cite{SW} that for $n \geq 1$, $p,r>1$, 
\begin{equation}\label{eq-SW-L}
0<\lambda <n
\end{equation}
and
\begin{equation}\label{eq-SW-AB}
\alpha < n(p-1)/p, \quad \beta < n (r-1)/r
\end{equation}
with 
\begin{equation}\label{eq-SW-A+B}
\alpha + \beta \geq 0,
\end{equation}
the following doubly weighted inequality holds
\begin{equation}\label{eq-SW}
\Big |\iint_{\R^n \times \R^n} \frac{f(x)g(y)}{|x|^\alpha |x-y|^\lambda |y|^\beta} dx dy\Big| 
\lesssim \| f \| _{L^p(\R^n )} \| g\| _{L^r(\R^n)}
\end{equation}
for any $f\in L^p(\R^n)$ and any $g\in L^r(\R^n)$ under the following conditions
\begin{equation}\label{eq-SW-Key}
1/p+1/r \geq 1
\end{equation}
and the balance condition
\begin{equation}\label{eq-SW-BC}
1 / p + 1/ r + (\lambda +\alpha+\beta)/n =2.
\end{equation}
In the literature, the inequality \eqref{eq-SW} is widely known as the Stein--Weiss inequality. Apparently, the SW inequality \eqref{eq-SW} becomes the classical HLS inequality \eqref{eq-HLS} if $\alpha = \beta =0$. In the original work of Stein and Weiss \cite{SW}, the condition \eqref{eq-SW-Key} is replaced by the following condition
\[
p \leq q := \Big(1 - \frac 1 r \Big)^{-1}.
\]
However, it is easy to verify that the above condition $p \leq q$ is actually equivalent to \eqref{eq-SW-Key}. Compared to the situation of \eqref{eq-HLS}, in the case of \eqref{eq-SW} there are three extra conditions for $\alpha$ and $\beta$: \eqref{eq-SW-AB}, \eqref{eq-SW-A+B}, and \eqref{eq-SW-Key}. It is clear that these three conditions are automatically satisfied if $\alpha=\beta=0$. 

Our motivation of writing this note comes from the fact that there was no further information on the conditions \eqref{eq-SW-AB}, \eqref{eq-SW-A+B}, and \eqref{eq-SW-Key} in the work of Stein and Weiss. It is worth noting that in the case $n=1$, Hardy and Littewood showed that the three conditions \eqref{eq-SW-AB}, \eqref{eq-SW-A+B}, and \eqref{eq-SW-Key} are \textit{necessary}. From this difference, a natural question is that \textit{whether these conditions are also necessary}. 

In the existing literature, these three conditions are often assumed, without justification, in many works such as \cite{HLZ12} for the case of the Heisenberg group, \cite{Dou2016} for the case of the upper half space, \cite{CLLT} for the case with the fractional Poisson kernel, and \cite{KRS} for the case of homogeneous Lie groups.

Our first observation in this note is the following which indicates the necessity of all conditions appearing in the SW inequality \eqref{eq-SW}.

\begin{theorem}\label{thm-MAIN}
All the conditions \eqref{eq-SW-L}, \eqref{eq-SW-AB}, \eqref{eq-SW-A+B}, \eqref{eq-SW-Key}, and \eqref{eq-SW-BC} are necessary for the Stein--Weiss inequality \eqref{eq-SW}.
\end{theorem}

We prove Theorem \ref{thm-MAIN} in section \ref{sec-Proof} by constructing counter-examples whenever one of these conditions is violated. Clearly, our result now concludes that the SW inequality \eqref{eq-SW} in its full generality. A careful examination of the necessary conditions for \eqref{eq-SW} reveals that the quantity $\alpha+\beta$ cannot be too positive. In fact, the size of $\alpha +\beta$ is implicitly encoded in \eqref{eq-SW-BC} because it follows from \eqref{eq-SW-BC} that
\[
\alpha+\beta \leq n -\lambda,
\]
which cannot be obtained directly from \eqref{eq-SW-AB}. Furthermore, it is clear that we have the following
\[
\text{\eqref{eq-SW-A+B}, \; \eqref{eq-SW-Key}, \; and \; \eqref{eq-SW-BC}} \quad \longrightarrow \quad \text{\eqref{eq-SW-L}}.
\] 
Therefore, to conclude the necessity of all conditions for the SW inequality \eqref{eq-SW}, we must not use the other condition as a means.

To further illustrate the idea of constructing counter-examples when proving Theorem \ref{thm-MAIN}, let us consider the doubly weighted Hardy--Littlewood--Sobolev inequality on the upper half space $\R_+^n$ which is defined as follows
$$\R_+^n = \{ (x',x_n) : x_n > 0\}.$$ 

In 2016, Dou proved the so-called the doubly weighted HLS inequality on $\R^{n-1} \times \R_+^n$; see \cite{Dou2016} and \cite{DZ-IMRN}. To be precise, given $n \geq 2$, $p,r>1$,
\begin{equation}\label{eq-SW-L-u}
0<\lambda <n-1 
\end{equation}
and
\begin{equation}\label{eq-SW-AB-u}
\alpha < (n-1)(p-1)/p, \quad \beta < n (r-1)/r
\end{equation}
with 
\begin{equation}\label{eq-SW-A+B-u}
\alpha + \beta \geq 0,
\end{equation} 
the following inequality holds
\begin{equation}\label{eq-SW-u}
\Big |\iint_{\R_+^n \times \R^{n-1}} \frac{f(x)g(y)}{|x|^\alpha |x-y|^\lambda |y|^\beta} dx dy\Big| 
\lesssim \| f \| _{L^p(\R^{n-1} )} \| g\| _{L^r(\R_+^n)}
\end{equation}
for any $f\in L^p(\R^{n-1})$ and any $g\in L^r(\R_+^n)$ under the following conditions
\begin{equation}\label{eq-SW-Key-u}
 1 / p + 1 / r \geq 1
\end{equation}
and the balance condition
\begin{equation}\label{eq-SW-BC-u}
\frac{n-1}n \frac 1 p + \frac 1 r + \frac{\lambda +\alpha+\beta + 1}n =2.
\end{equation}
Here we identify $\R^{n-1}$ as the boundary of $\R_+^n$ and the `distance' $|x-y|$ is understood as follows 
$$|x-y| = \sqrt{|x-y'|^2 + y_n^2}.$$ 
Compared to the classical SW inequality \eqref{eq-SW} on $\R^n \times \R^n$, the two conditions \eqref{eq-SW-A+B-u} and \eqref{eq-SW-Key-u} remains the same as those of \eqref{eq-SW}, only the three conditions \eqref{eq-SW-L-u}, \eqref{eq-SW-AB-u}, and \eqref{eq-SW-BC-u} change. However, we \textit{no longer} have the following
\[
\text{\eqref{eq-SW-A+B-u},\; \eqref{eq-SW-Key-u}, \; and \; \eqref{eq-SW-BC-u}} \quad \longrightarrow \quad \text{\eqref{eq-SW-L-u}}
\] 
as in the previous case. This suggests that one of the three conditions \eqref{eq-SW-L-u}, \eqref{eq-SW-AB-u}, and \eqref{eq-SW-BC-u} should not be necessary. Based on this intriguing observation, we are forced to further investigate this by seek for the necessity of the above conditions. The following result is what we obtain.

\begin{theorem}\label{thm-MAIN-u}
All the conditions \eqref{eq-SW-L-u} replaced by $0<\lambda<n-1/r$, \eqref{eq-SW-AB-u}, \eqref{eq-SW-A+B-u}, \eqref{eq-SW-Key-u}, and \eqref{eq-SW-BC-u} are necessary for the Stein--Weiss inequality \eqref{eq-SW-u}.
\end{theorem}

We prove Theorem \ref{thm-MAIN-u} in section \ref{sec-Proof-u} below. The new threshold $n-1/r$ for $\lambda$ appears by means of the identity \eqref{eq-SW-BC-u-n1} below, which as far as we know has not been observed before. As can be seen from Theorem \ref{thm-MAIN-u}, we cannot say anything if 
$$n-1 \leq \lambda \leq n-1/r.$$ 
In fact, one could expect by seeing the condition \eqref{eq-SW-L-u} used in \cite{Dou2016} that one still needs $\lambda < n-1$. 
From this fundamental difference, it is natural to ask whether or not the SW inequality \eqref{eq-SW-u} on $\R_+^n$ is \textit{still valid} for the following new range of the parameter $\lambda$
$$n-1 \leq \lambda < n-1/r.$$ 
In other words, if we fix $0<\lambda<n$, this is equivalent to saying that the inequality \eqref{eq-SW-u} remains valid so long as
\[
r > \max \Big\{ 1, \frac 1{n-\lambda} \Big\}.
\]
It turns out that this is indeed the case. In the paper \cite{NNN}, served as the second in the set, we address this issue. To be more precise, we prove in \cite{NNN} the inequality \eqref{eq-SW-m} below, namely we are interested in the validity of the following inequality
\begin{equation}\label{eq-SW-m}
\Big |\iint_{\R^n \times \R^{n-k}} \frac{f(x)g(y)}{|x|^\alpha |x-y|^\lambda |y|^\beta} dx dy\Big| 
\lesssim \| f \| _{L^p(\R^{n-k} )} \| g\| _{L^r(\R^n)}
\end{equation}
for any $f\in L^p(\R^{n-k})$ and any $g\in L^r(\R^n)$. Here we assume $n > k \geq 0$. Clearly, the above inequality \eqref{eq-SW-m} includes \eqref{eq-SW-u} as a special case if we choose $k=1$. 
As far as we know, the inequality \eqref{eq-SW-m} is not yet studied. Although the two inequalities \eqref{eq-SW-u} and \eqref{eq-SW-m} look rather similar, the motivation of working on the inequality \eqref{eq-SW-m} is not just to create an artificial one, but it stems from an intriguing connection between isoperimetric inequalities and HLS inequalities with Poisson-type kernel on $\R^{n-1} \times \R_+^n$; see \cite{HWY08} for the intriguing connection between these inequalities. It is shown in \cite{NNN} that the inequality \eqref{eq-SW-m} provides more information on HLS inequalities with Poisson-type kernel on $\R^{n-1} \times \R_+^n$ and this is the key observation to work on \eqref{eq-SW-m}. 

Toward the validity of the inequality \eqref{eq-SW-m}, we aim to address the necessity for the inequality \eqref{eq-SW-m}, leading to the content of the next result of the present paper. Inspired by the set of conditions for \eqref{eq-SW-u}, our set of assumptions for \eqref{eq-SW-m} is as follows. First we let $p,r>1$,
\begin{equation}\label{eq-SW-L-m}
0<\lambda <n-k/r ,
\end{equation}
and
\begin{equation}\label{eq-SW-AB-m}
\alpha < (n-k)(p-1)/p, \quad \beta < n (r-1)/r,
\end{equation}
with
\begin{equation}\label{eq-SW-A+B-m}
\alpha + \beta \geq 0.
\end{equation} 
We also assume the following conditions
\begin{equation}\label{eq-SW-Key-m}
 1 / p + 1 / r \geq 1
\end{equation}
and
\begin{equation}\label{eq-SW-BC-m}
\frac{n-k}n \frac 1 p + \frac 1 r + \frac{\lambda +\alpha+\beta + k}n =2.
\end{equation}
Now we identify necessary conditions for the inequality in the same fashion of the previous two inequalities \eqref{eq-SW} and \eqref{eq-SW-u}. We shall prove the following necessity result.

\begin{theorem}\label{thm-MAIN-m}
All the conditions \eqref{eq-SW-L-m}, \eqref{eq-SW-AB-m}, \eqref{eq-SW-A+B-m}, \eqref{eq-SW-Key-m}, and \eqref{eq-SW-BC-m} are necessary for the validity of the Stein--Weiss inequality \eqref{eq-SW-m}.
\end{theorem}

We prove Theorem \ref{thm-MAIN-m} in section \ref{sec-Proof-m} below. As in the case of the upper half space, we obviously do \textit{not} have the following
\[
\text{\eqref{eq-SW-A+B-m},\; \eqref{eq-SW-Key-m}, \; and \; \eqref{eq-SW-BC-m}} \quad \longrightarrow \quad \text{\eqref{eq-SW-L-m}}.
\] 
In \cite{NNN}, we investigate the validity of \eqref{eq-SW-m} under the above necessary conditions. To be more precise, instead of proving \eqref{eq-SW-m} in its form, we simply consider \eqref{eq-SW-m} with $\alpha=0$ and the weight $|y|^{-\beta}$ replaced by $|y''|^{-\beta}$, namely
\begin{equation}\label{eq-SW-m-simpler}
\Big |\iint_{\R^n \times \R^{n-k}} \frac{f(x)g(y)}{ |x-y|^\lambda |y''|^\beta} dx dy\Big| 
\lesssim \| f \| _{L^p(\R^{n-k} )} \| g\| _{L^r(\R^n)}
\end{equation}
with $y=(y',y'') \in \R^{n-k} \times \R^k$. Clearly, a special case of the inequality \eqref{eq-SW-m-simpler} is the following HLS inequality on $\R^{n-k} \times \R^n$ with $n > k \geq 0$
\begin{equation}\label{eq-HLS-m}
\Big |\iint_{\R^n \times \R^{n-k}} \frac{f(x)g(y)}{ |x-y|^\lambda } dx dy\Big| 
\lesssim \| f \| _{L^p(\R^{n-k} )} \| g\| _{L^r(\R^n)}.
\end{equation}
Surprisingly, with the new technique introduced in \cite{NNN}, the validity of \eqref{eq-SW-m} can be quickly derived from \eqref{eq-HLS-m}. This provides a very short proof of the SW inequality \eqref{eq-SW-u} on the upper half space $\R_+^n$. For the classical approach of Stein--Weiss applied to the case of the upper half space $\R_+^n$ with the limitation $0<\lambda<n-1$, we refer the reader to the work \cite{Dou2016}. We expect that the classical approach of Stein--Weiss can also be applied to the case of \eqref{eq-SW-m} giving another the proof for $0<\lambda<n-k$, but we are not so sure if this approach can cover the range $n-k \leq \lambda < n-k/r$. We also expect that the approach in \cite{SDY2014}, which is based on interpolation of Lorentz spaces, can also be applied in this case. Hence, we might have another proof of \eqref{eq-SW-m}.

Before closing this section, we have three comments. First, we can easily show that all the conditions for the Stein--Weiss inequality
\[
\Big |\iint_{\mathbb G \times \mathbb G} \frac{f(x)g(y)}{|x|^\alpha |y^{-1}x|^\lambda |y|^\beta} dx dy\Big| 
\lesssim \| f \| _{L^p (\mathbb G) } \| g\| _{L^r (\mathbb G) }
\]
on the product of the homogeneous Lie group $\mathbb G$ are also necessary; see \cite{KRS} for the precise statement of the inequality. Since the argument in the proof of Theorem \ref{thm-MAIN} can be used for the context of homogeneous Lie groups directly without any difficulty, we omit the details.

For the second comment, over the last few years, there have been intensive studies on reverse cases of the classical HLS inequality and well as the classical SW inequality. For example, the authors in \cite{CLLT-TAMS} consider the reverse SW inequality on $\R^n \times \R^n$, following the case on $\R^{n-1} \times \R_+^n$ in \cite{CLT-ANS}. These inequalities can formally be written as follows
\begin{equation}\label{eq-SW-r}
\Big |\iint \frac{f(x)g(y)}{|x|^\alpha |x-y|^\lambda |y|^\beta} dx dy\Big| 
\gtrsim \| f \| _{L^p } \| g\| _{L^r }
\end{equation}
with $\lambda <0$ and $p, r \in (0,1)$. Accordingly, the above inequality also enjoys the same balance condition as in \eqref{eq-SW-BC} for the case $\R^n \times \R^n$ and \eqref{eq-SW-BC-u} for case $\R^{n-1} \times \R_+^n$. However, unlike the inequalities \eqref{eq-SW} and \eqref{eq-SW-u}, the reverse case \eqref{eq-SW-r} considered in the above mentioned papers does not exhibit similar `necessary' conditions except a suitable balance condition and the natural assumption $\lambda <0$. Using the balance condition and due to the fact that $\lambda<0$, there are foreseen conditions
\[
\alpha + \beta > n [ (p-1)/p + (r-1)/r ]
\]
for the case $\R^n \times \R^n$ and
\[
\alpha + \beta > n (p-1)/p + (n-1)(r-1)/r
\]
for the case $\R^{n-1} \times \R_+^n$. Consequently, to verify the validity of \eqref{eq-SW-r} it is natural to assume that at least one of the following two conditions
\[
\alpha >n(p-1)/p, \quad \beta > n(r-1)/r
\]
for the case $\R^n \times \R^n$ and
\[
\alpha > n(p-1)/p , \quad \beta > (n-1)(r-1)/r 
\]
for the case $\R^{n-1} \times \R_+^n$ occurs. 

However, we can say more about the two sets of conditions above. Following the argument proving the necessity of the conditions \eqref{eq-SW-AB} and \eqref{eq-SW-AB-u} established in Theorems \ref{thm-MAIN} and \ref{thm-MAIN-u} above, it is quite clear that if at least one of the two conditions for $\alpha$ and $\beta$ for each case does not not occur, the left hand side of \eqref{eq-SW-r} diverges, leading to the triviality of \eqref{eq-SW-r}. This sheds light on why the authors in \cite[Theorem 1]{CLLT-TAMS} can assume that the following conditions
\[
0 \geq \alpha > n(p-1)/p, \quad 0 \geq \beta > n(r-1)/r.
\]
However, we suspect that the negativity of both $\alpha$ and $\beta$ seems to be not necessary.

As the final remark before closing this section, we would like to note that during the preparation of this work, it has come to our attention that the necessity of all conditions for the SW inequality \eqref{eq-SW} was already studied in \cite{SDY2014}. Nevertheless and very fortunate, the argument and examples demonstrated in \eqref{eq-SW} are essentially different from ours. We also expect that our finding in this note could be useful and applicable for other cases. 


This note is organized as follows:

\renewcommand{\baselinestretch}{0.75}\normalsize
\tableofcontents
\renewcommand{\baselinestretch}{1.0}\normalsize

From now on, by $B_r^\ell$ we mean the open ball in $\R^\ell$ centered at the origin with radius $r$. Since we are also interested in the upper half space $\R_+^n$, for simplicity by $B_r^n$ we also mean $\R_+^n \cap B_r^n$ if no confusion occurs.


\section{Necessity of conditions for (\ref{eq-SW}): proof of Theorem \ref{thm-MAIN}}
\label{sec-Proof}

In this section, we give a proof of Theorem \ref{thm-MAIN}. As mentioned earlier, the result is already known, thanks to \cite{SDY2014}. Hence, we aim to provide new examples to support the finding, especially to be able to consider the case of the upper half space $\R^{n-1} \times \R_+^n$ as well as $\R^{n-k} \times \R^n$. 

Compared to the examples constructed in \cite{SDY2014}, our main contribution in this section is a new example constructed in subsection \ref{subsec-A+B} below, which is inspired by a similar construction due to Hardy and Littewood.


\subsection{The necessity of (\ref{eq-SW-L}), (\ref{eq-SW-AB}), and (\ref{eq-SW-BC})}
\label{subsec-L}


\subsubsection{The necessity of $0<\lambda<n$}

We start with the necessity of \eqref{eq-SW-L}, which is quite obvious and well-known. Indeed, first we rule out the case $\lambda \leq 0$. Clearly for any each fixed $x$, near the infinity the quantity $|x-y|^{-\lambda}$ is as large as we want if $\lambda < 0$. Hence \eqref{eq-SW} cannot hold if $\lambda < 0$. The case $\lambda = 0$ cannot occur because
\[
\iint_{\R^n \times \R^n} \frac{dxdy}{|x|^\alpha |y|^\beta} > \iint_{\{ 1 \leq |x| \leq 2\} \times \R^n} \frac{dxdy}{|x|^\alpha |y|^\beta} = +\infty
\]
regardless of $\beta$. It is worth noticing that the case $\lambda = 0$ is often called the limiting case and in this scenario, we expect that there is a log-SW inequality in the same fashion of the log-HLS inequality; see \cite[Corollary 5.3]{DZ-IMRN}. Clearly, by seeing the necessity of the conditions \eqref{eq-SW-A+B} and \eqref{eq-SW-Key}, we must have $\lambda \leq n$. However, we still need to rule out the case $\lambda = n$. To rule out the case $\lambda \geq n$, we simply choose
\[
f = g = \chi_{B_6 \setminus B_1}.
\]
and show that the two terms $|x|^{-\alpha}$ and $|y|^{-\beta}$ are negligible. Then, by Tonelli's theorem, it is easy to see that
\begin{align*}
+\infty 
&> \iint_{\R^n \times \R^n} \frac{\chi_{B_6 \setminus B_1} (x) \chi_{B_6 \setminus B_1} (y)}{|x|^\alpha |x-y|^\lambda |y|^\beta} dx dy \\
&\gtrsim \int_{B_4 \setminus B_2} \Big( \int_{B_{|y|/2}(y)} \frac{dx}{|x-y|^\lambda} \Big) dy \\
\overset{|y| \geq 2} &\geq \int_{B_4 \setminus B_2} \Big( \int_{B_1 (y)} \frac{dx}{|x-y|^\lambda} \Big) dy 
= +\infty.
\end{align*}
This completes the proof of the necessity of $\lambda < n$.

\begin{remark}\label{rmk-SW-L}
In the above argument, if we use the dual version of \eqref{eq-SW-u}, namely
\begin{equation}\label{eq-Duality-f}
\| f \|_{L^p(\R^n)}^q \geq \int_{\R^n} \Big( \int_{\R^n} \frac{f(x)}{ |x|^\alpha |x-y|^\lambda |y|^\beta } dx \Big)^q dy
\end{equation}
with $q = r/(r-1)$, then the above argument still works. Indeed, we easily obtain
\begin{align*}
\int_{\R^n} \Big( \int_{\R^n} \frac{ \chi_{B_6 \setminus B_1} (x)}{ |x|^\alpha |x-y|^\lambda |y|^\beta } dx \Big)^q dy 
&\gtrsim \int_{B_4 \setminus B_2} \Big( \int_{B_{|y|/2}(y)} \frac{dx}{|x|^\alpha |x-y|^\lambda} \Big)^q dy \\
&\gtrsim \int_{B_4 \setminus B_2} \Big( \int_{B_1 (y)} \frac{dx}{|x-y|^\lambda} \Big)^q dy 
= +\infty
\end{align*}
if $\lambda \geq n$.
\end{remark}



\subsubsection{The necessity of $\beta < n(r-1)/r$}
\label{subsubsec-beta<}

First, we prove the necessity of the condition $\beta < n(r-1)/r$. By way of contradiction, we assume 
$$\beta \geq n(r-1)/r=n/q.$$ 
We shall obtain contradiction from \eqref{eq-Duality-f}. Then we choose $f = \chi_{B_3 \setminus B_2}$ and show that the two terms $|x|^{-\alpha}$ and $|x-y|^{-\lambda}$ in \eqref{eq-Duality-f} are negligible. Indeed, we can estimate
\begin{align*}
\int_{\R^n} \Big( \int_{\R^n} \frac{\chi_{B_3 \setminus B_2 }(x)}{ |x|^\alpha |x-y|^\lambda |y|^\beta } dx \Big)^q dy
&\geq \int_{B_1 } 
\Big[ \int_{B_3 \setminus B_2} \frac {dx }{ |x|^\alpha |x-y|^\lambda} \Big]^q 
\frac{dy}{|y|^{\beta q}} 
= +\infty.
\end{align*}
This violates \eqref{eq-Duality-f}, hence concluding the necessity of $\beta < n(r-1)/r$. Notice that in the above estimate, we also need the uniformly lower bound of $\int_{B_3 \setminus B_2 } |x|^{-\alpha} |x-y|^{-\lambda}dx$ but this is quite obvious because $2 \leq |x| \leq 3$ and $$1 \leq |x|-|y| \leq |x-y| \leq |x| + |y| \leq 4.$$


\subsubsection{The necessity of $\alpha < n(p-1)/p$}
\label{subsubsec-alpha<}

Now we establish the necessity of the condition $\alpha < n(p-1)/p$. In this case, we use the following equivalent form of \eqref{eq-SW}
\begin{equation}\label{eq-Duality-g}
\| g \|_{L^p(\R^n)}^q \geq \int_{\R^n} \Big( \int_{\R^n} \frac{g(y)}{ |x|^\alpha |x-y|^\lambda |y|^\beta } dy \Big)^q dx
\end{equation}
with $q = p/(p-1)$. By contradiction, we assume 
$$\alpha \geq n(p-1)/p=n/q.$$
Then we choose $g = \chi_{B_3 \setminus B_2}$ and show that the two terms $|x-y|^{-\lambda}$ and $|y|^{-\beta}$ are negligible. Indeed, we can estimate
\begin{align*}
\int_{\R^n} \Big( \int_{\R^n} \frac{\chi_{B_3 \setminus B_2}(y) }{ |x|^\alpha |x-y|^\lambda |y|^\beta } dy \Big)^q dx
&\geq \int_{B_1 } 
\Big[ \int_{B_3 \setminus B_2} \frac {dy }{ |x-y|^\lambda |y|^\beta } \Big]^q 
\frac{dx}{|x|^{\alpha q}} = +\infty.
\end{align*}
This violates \eqref{eq-Duality-g}, hence concluding the necessity of $\alpha < n(p-1)/p$. 


\subsubsection{The necessity of $1/p + 1/r + (\lambda + \alpha + \beta)/n =2$}
\label{subsubsec-=}

This part is well-known due to the scaling invariance of the inequality. Since its argument is short, we include its proof for completeness. Fix two positive functions $f$ and $g$ such that $\|f\|_{L^p(\R )}<+\infty$ and $\|g\|_{L^r(\R^n)}<+\infty$. Denote
\begin{equation}\label{eq-FeGe}
f_\epsilon (x) = \epsilon^{-\frac np} f \big( \frac x \epsilon \big),
\quad
g_\epsilon (y) = \epsilon^{-\frac nr} g \big( \frac y \epsilon \big).
\end{equation}
Obviously,
\[
\|f_\epsilon\|_{L^p(\R )} = \|f\|_{L^p(\R )}<+\infty, \quad
\|g_\epsilon\|_{L^r(\R^n)} = \|g\|_{L^r(\R^n)}<+\infty.
\]
A simple change of variables gives
\begin{align*}
\iint_{\R^n \times \R^n } \frac{f_\epsilon (x) g_\epsilon (y) dx dy}{|x|^\alpha |x-y|^\lambda |y|^\beta } 
&=\iint_{\R^n \times \R^n } \frac{f (\frac x \epsilon ) g (\frac y \epsilon ) dx dy }{\epsilon^{\frac np + \frac nr + \lambda + \alpha + \beta} |\frac x\epsilon|^\alpha |\frac x\epsilon - \frac y \epsilon |^\lambda |\frac{y}\epsilon|^\beta } \\
&=\big( \frac 1 \epsilon \big)^{( \frac 1p + \frac 1r + \frac{\lambda + \alpha + \beta}n - 2 ) n}\iint_{\R^n \times \R^n } \frac{f (x ) g (y ) dx dy}{|x|^\alpha |x -y |^\lambda |y |^\beta } .
\end{align*}
Hence, if the balance condition \eqref{eq-SW-BC} does not hold, we simply send $\epsilon$ to zero or to infinity to get a contradiction.

\begin{remark}
J. Dou kindly informs us that the necessity of \eqref{eq-SW-BC} can also be derived from a suitable Pohozaev-type identity. We leave the detail for the interested readers.
\end{remark}


\subsection{The necessity of (\ref{eq-SW-A+B})}
\label{subsec-A+B}

To verify the necessity of $\alpha + \beta \geq 0$ for the SW inequality on $\R^n$, we construct an example of non-negative functions $f$ and $g$ in such a way that $\|f\|_{L^p(\R^n)} <+\infty$, $\|g\|_{L^r(\R^n)} <+\infty$, but
\[
\iint_{\R^n \times \R^n} \frac{f(x)g(y)}{|x|^\alpha |x-y|^\lambda |y|^\beta} dx dy = +\infty.
\] 
Our example is inspired by a similar example in 1-dimensional case due to Hardy and Littewood; see \cite[page 580]{hl1928}. Roughly speaking, we construct an example in such a way that the following two properties are satisfied:
\begin{itemize}
 \item $|x-y|^{-\lambda}$ is bounded from below away from zero and
 \item $|x|$ and $|y|$ are compatible (near infinity).
\end{itemize}
For a general point $x \in \R^n$, we often write $x=(x',x_n) \in \R^{n-1} \times \R$. Now fix any small $\epsilon > 0$ and define
\[
f(x)=
\left\{
\begin{aligned}
& (\log_2 x_n)^{-\frac 1p - \epsilon} &&\text{if } \; 2^m \leq x_n \leq 2^m+1, \; m \geq 1,\\
& & & \text{and } \; |x'| \leq 1,\\
&0 &&\text{otherwise},
\end{aligned}
\right.
\]
and
\[
g(y)=
\left\{
\begin{aligned}
& (\log_2 y_n)^{-\frac 1r - \epsilon} &&\text{if } \; 2^m \leq y_n \leq 2^m+1, \; m \geq 1, \\
& & & \text{and } \; |y'| \leq 1,\\
&0 &&\text{otherwise}.
\end{aligned}
\right.
\]
Intuitively, the support of $f$ and $g$ consists of \textit{disjoint cylinders} concentrated along the half axis $x_n$. The condition $|x'| \leq 1$ is used to make the functions $f$ and $g$ essentially depending only on the last coordinate. Observe that
\begin{align*}
\int_{\R^n} f(x)^p dx &=\sum_{m \geq 1} \int_{2^m}^{2^m+1}\Big( \int_{|x'| \leq 1} \frac 1{(\log_2 x_n)^{1+ p\epsilon}} dx' \Big) dx_n\\
&\leq |B_1^{n-1} | \sum_{m \geq 1} \frac{1}{m^{1+ p\epsilon}}
<+\infty.
\end{align*}
Hence $f \in L^p(\R^n)$. Similarly, we also obtain $g \in L^r(\R^n)$. Now we estimate the double integral as follows
\begin{equation}\label{eq-Double-1}
\begin{aligned}
\iint_{\R^n \times \R^n} & \frac {f(x)g(y)}{|x|^\alpha |x-y|^\lambda |y|^\beta} dxdy \\
&\geq\sum_{m \geq 1} \iint_{B_1^{n-1} \times [2^m, 2^m+1]} \Big( \iint_{B_1^{n-1} \times [2^m, 2^m+1]} \frac {g(y) dy}{ |x-y|^\lambda |y|^\beta} \Big) \frac{f(x) dx}{|x|^\alpha}.
\end{aligned}
\end{equation}
For each $m \geq 1$, let $x$ and $y$ be such that
\begin{equation}\label{eq-Strip}
x,y \in \overline{B_1^{n-1}} \times [2^m, 2^m+1].
\end{equation}
We need to estimate $|x|^{-\alpha}$, $|y|^{-\beta}$, and $|x-y|^{-\lambda}$ from below. Apparently, we easily bound $|x-y|$ from the above as follows
\[
|x-y|^2 =|x'-y'|^2 + |x_n-y_n|^2 \leq 2(|x'|^2+|y'|^2) + 1 \leq 5.
\]
This together with $\lambda > 0$ allows us to estimate
\[
|x-y|^{-\lambda} \geq 3^{-\lambda}
\]
provided \eqref{eq-Strip} holds for each $m \geq 1$. Since $|x'| \leq 1$ and $x_n \geq 2$, it is easy to verify that
\[
x_n \leq \sqrt{|x'|^2+x_n^2} \leq 2x_n.
\]
Similarly, we get
\[
y_n \leq \sqrt{|y'|^2+y_n^2} \leq 2y_n.
\]
Hence, regardless of the sign of $\alpha$ and $\beta$, we always have
\[
|x|^{-\alpha} \gtrsim x_n^{-\alpha} \quad \text{and} \quad
|y|^{-\beta} \gtrsim y_n^{-\beta}.
\]
Putting the above estimates together, we can further estimate \eqref{eq-Double-1} as follows
\begin{equation}\label{eq-Double-2}
\begin{aligned}
\iint_{\R^n \times \R^n} & \frac {f(x)g(y)}{|x|^\alpha |x-y|^\lambda |y|^\beta} dxdy\\
& \gtrsim \sum_{m \geq 1} \Big( \iint_{B_1^{n-1} \times [2^m, 2^m+1]} \frac{f(x) dx}{|x|^\alpha} \Big) \Big( \iint_{B_1^{n-1} \times [2^m, 2^m+1]} \frac {g(y) dy}{ |y|^\beta} \Big) \\
& = |B_1^{n-1}|^2 \sum_{m \geq 1} \Big( \int_{2^m}^{2^m+1} \frac {dx_n}{x_n^\alpha (\log_2 x_n)^{\frac 1p +\epsilon}} \Big) \Big( \int_{2^m}^{2^m+1} \frac {dy_n}{y_n^\beta (\log_2 y_n)^{\frac 1r +\epsilon}} \Big) \\
& \gtrsim\sum_{m \geq 1} \frac 1{(m+1)^{\frac 1p +\frac 1r +2\epsilon}} 
\Big( \int_{2^m}^{2^m+1} \frac {dx_n}{x_n^\alpha } \Big) \Big( \int_{2^m}^{2^m+1} \frac {dy_n}{y_n^\beta } \Big).
\end{aligned}
\end{equation}
Notice that for $2^m \leq x_n \leq 2^m+1$ we always have
\[
\frac 1{x_n^\alpha} \geq 
\left\{
\begin{aligned}
& \frac 1{2^\alpha} \frac 1{(2^m)^\alpha} &&\text{if } \; \alpha \geq 0,\\
& \frac 1{(2^m)^\alpha} &&\text{if } \; \alpha < 0.
\end{aligned}
\right.
\]
Similarly, we also obtain $y_n^{-\beta} \gtrsim (2^m)^{-\beta}$ provided $2^m \leq y_n \leq 2^m+1$. Hence
\[
\Big( \int_{2^m}^{2^m+1} \frac {dx_n}{x_n^\alpha } \Big) \Big( \int_{2^m}^{2^m+1} \frac {dy_n}{y_n^\alpha } \Big)
\geq \frac 1{(2^m)^{\alpha + \beta}}.
\]
Thus, recalling \eqref{eq-Double-2} we arrive at
\begin{align*}
\iint_{\R^n \times \R^n} & \frac {f(x)g(y)}{|x|^\alpha |x-y|^\lambda |y|^\beta} dxdy
\gtrsim \sum_{m \geq 1} \frac 1{(2^m)^{\alpha + \beta}} \frac 1{(m+1)^{\frac 1p +\frac 1r + 2\epsilon}} .
\end{align*}
Obviously, one can easily see that the above sum diverges if $\alpha +\beta < 0$. Hence, the necessity of \eqref{eq-SW-AB} is proved.


\subsection{The necessity of (\ref{eq-SW-Key})}
\label{subsec-P+R}

Now we establish the necessity of the condition \eqref{eq-SW-Key}. Indeed, suppose $1/p+1/r<1$. Then
\[
q := \Big(1- \frac 1r \Big)^{-1} < p.
\] 
(Note that this is also equivalent to $\lambda+\alpha+\beta>n$.) In view of \eqref{eq-Duality-f}, it suffices to show that the left hand side of \eqref{eq-Duality-f} becomes infinity for suitable $f \in L^p(\R^n)$. By modifying the counter-example constructed in subsection \ref{subsec-A+B} above, we now choose 
\[
f(x)=
\left\{
\begin{aligned}
& |x|^{-\frac np} (\log |x|)^{-\frac 1 q} & & \text{if } |x| \geq 2,\\
& 0 & & \text{otherwise}.
\end{aligned}
\right.
\]
It is not hard to see that
\[
\int_{\R^n} f(x)^p dx = |\mathbb S^{n-1}| \int_2^{+\infty} \frac{d\rho}{\rho (\log \rho)^{p/q}} < +\infty,
\]
thanks to $p>q$. Hence, $f \in L^p(\R^n)$. To obtain a contradiction, we estimate the double integral near the region $\{|x-y| = 0\}$; see \cite[page 580]{hl1928}. By \eqref{eq-SW} we have
\begin{align*}
+\infty > \|f\|_{L^p(\R^n)}^q &> \int_{\R^n } \Big( \int_{\R^n} \frac{ f(x) }{|x|^\alpha |x-y|^\lambda |y|^\beta} dx \Big)^q dy\\
&\geq
\int_{|y| \geq 4} \Big( \int_{|y|/2 \leq |x| \leq 2|y|} \frac{dx}{|x|^{\alpha + \frac np} (\log |x|)^{ \frac 1 q} |x-y|^\lambda } \Big)^q \frac {dy}{|y|^{\beta q}}.
\end{align*}
As $\lambda > 0$ and $|x| \leq 2|y|$, we have $|x-y|^{\lambda} \leq 3^\lambda |y|^\lambda$, which then yields
\[
|x|^{\alpha + \frac np} |x-y|^{\lambda} \leq 2^{\alpha + \frac np} 3^\lambda |y|^{\lambda + \alpha + \frac np}
\]
for any $|x|/2 \leq |y| \leq 2|x|$. (We need the extra condition $|y| \geq |x|/2$ because $\alpha + n/p$ could be negative.) Hence, together with 
\[
0< \log |x| \leq \log (2|y|) \leq \log 2 + \log |y| \leq 2\log |y|,
\] 
we arrive at
\[
\frac 1{|x|^{\alpha + \frac np} (\log |x|)^{ \frac 1 q} |x-y|^\lambda }
\gtrsim \frac 1{|y|^{ \lambda + \alpha + \frac np} (\log |y|)^{ \frac 1 q} }.
\]
Keep in mind that $q>1$. Thus,
\begin{align*}
+\infty > \|f\|_{L^p(\R^n)}^q
&\gtrsim
\int_{|y| \geq 4} \Big( \int_{|y|/2 \leq |x| \leq 2|y|} dx \Big)^q \frac{ dy}{|y|^{(\lambda+\alpha +\beta+ \frac np ) q } \log |y|}\\
& \gtrsim \int_{|y| \geq 4} \frac{ dy}{|y|^n \log |y|}
= +\infty,
\end{align*}
thanks to
\begin{align*}
\Big(\lambda + \alpha +\beta + \frac np \Big) q &= \Big(\frac{\lambda + \alpha +\beta}n + \frac 1p \Big) n q 
= \Big( 1 - \frac 1r + 1 \Big) n q = n + nq.
\end{align*}
In the above estimate, we have used the fact that $\int_{|y|/2 \leq |x| \leq 2|y|} dx \gtrsim |y|^n$, thanks to $|y| \geq 4$. Hence, we obtain a contradiction.


\section{Necessity of conditions for (\ref{eq-SW-u}) on $\R_+^n \times \R^{n-1}$: proof of Theorem \ref{thm-MAIN-u}}
\label{sec-Proof-u}

In this section, we show that our argument performed in section \ref{sec-Proof} above can be applied to the case of the upper half space $\R_+^n$ with some changes. Although there are some similarities between the previous section and this section, for a number of reasons, the analysis in this section is considerably more difficult than that in section \ref{sec-Proof}.

For a general point $x \in \R^{n-1}$, we write $$x=(x',x_{n-1}) \in \R^{n-2} \times \R.$$ Hence for $y \in \R^n$ we can also write $$y = (y', y_{n-1}, y_n) \in \R^{n-2} \times \R \times \R.$$ 


\subsection{The necessity of $0<\lambda<n-1/r$}
\label{subsec-L-u}

Employing a similar argument as in subsection \ref{subsec-L}, we can easily rule out the case $\lambda \leq 0$. Now we rule out the case $\lambda \geq n-1/r$. Compared to the case of \eqref{eq-SW-L}, this case is non-trivial due to the presence of the factional number $n-1/r$. Our crucial observation can be described as follows. Still denote 
\[
q := \Big(1-\frac 1r\Big)^{-1},
\] 
the condition
\[
\lambda \geq n - 1/r \quad \text{is \textit{equivalent} to} \quad \lambda - 1/q \geq n-1.
\] 
It is important to note that the balance condition \eqref{eq-SW-BC-u} can be rewritten as
\begin{equation}\label{eq-SW-BC-u-n1}
\frac 1p + \frac 1r + \frac{\lambda + \alpha +\beta - \frac 1q}{n-1} = 2.
\end{equation}
Hence, the two conditions $\lambda - 1/q \geq n-1$ and $\alpha + \beta \geq 0$ basically say that we must have
\[
1/p + 1/r \leq 1,
\] 
which should not occur by seeing \eqref{eq-SW-Key-u}. However, we still need to rule out the case $1/p+1/r=1$, which corresponds to the case both
\[
\lambda - 1/q = n-1 \quad \text{and} \quad \alpha + \beta = 0
\] 
occur simultaneously by seeing \eqref{eq-SW-BC-u-n1}. This is a delicate issue, as we shall soon see. 

The difficulty comes from the presence of the fractional number $1/q$. Hence we need some new idea. It is worth noticing that our argument below works for all $\lambda \geq n-1/r$, not just $\lambda = n-1/r$, and it does not depend on other conditions for the inequality, as expected. 

To overcome this difficulty and in view of Remark \ref{rmk-SW-L}, we need to use the dual version of \eqref{eq-SW-u}. Clearly, for some non-negative function $f \in L^p(\R^{n-1})$ to be determined later, the inequality \eqref{eq-SW-BC-u} is equivalent to
\begin{align*}
+\infty > \| f \|_{L^p(\R^{n-1})}^q &\gtrsim \int_{\R_+^n} \Big( \int_{\R^{n-1}} \frac{f(x) }{ |x|^\alpha |x-y|^\lambda |y|^\beta } dx \Big)^q dy \\
&=\int_{\R^{n-1}} \int_0^{+\infty} \Big( \frac 1{ |y|^\beta } \int_{\R^{n-1}} \frac{f(x) }{ |x|^\alpha |x-y|^\lambda } dx \Big)^q dy_n dy',
\end{align*}
with $q=r/(r-1)$. Using Fubini's theorem and the non-decreasing of the integral $ \int_{B_{\rho}^{n-1} (y')} |x|^{-\alpha } f(x)dx$ in $\rho$, we easily get
\begin{align*}
\int_{\R^{n-1}} \frac{f(x) }{ |x|^\alpha |x-y|^\lambda } dx
&=\lambda \int_{\R^{n-1}} \Big( \int_{|x-y|}^{+\infty} \frac{d\rho}{\rho^{\lambda +1}} \Big) \frac{f(x) }{ |x|^\alpha } dx\\
&=\lambda \int_{y_n}^{+\infty} \Big( \int_{B_{\rho - y_n}^{n-1} (y')} \frac{f(x) }{ |x|^\alpha } dx \Big) \frac{d\rho}{\rho^{\lambda +1}}\\
&\gtrsim \frac 1{y_n^\lambda } \int_{B_{y_n}^{n-1} (y')} \frac{f(x) }{ |x|^\alpha } dx.
\end{align*}
Hence, we obtain
\begin{align*}
 \int_{\R^{n-1}} \int_0^{+\infty} & \Big( \frac 1{ |y|^\beta } \int_{\R^{n-1}} \frac{f(x) }{ |x|^\alpha |x-y|^\lambda } dx \Big)^q dy_n dy'\\
&\geq \int_{B_4^{n-1} \setminus B_2^{n-1}} \int_0^1 \Big( \frac 1{y_n^{ \lambda -1/q} |y|^\beta } \int_{B_{y_n}^{n-1} (y')} \frac{f(x) }{ |x|^\alpha } dx \Big)^q \frac{dy_n}{y_n} dy' .
\end{align*}
For the last line in the above computation, thanks to $|y'| \geq 2$ and $0\leq y_n \leq 1$, we know that 
\[
B_1^{n-1} \subset B_{y_n}^{n-1} (y') \subset B_5^{n-1}
\]
and that $2 \leq |y| \leq \sqrt 5$. Hence, if we choose $f = \chi_{B_6^{n-1}}$, then we can bound
\[
\frac 1{|y|^\beta } \int_{B_{y_n}^{n-1} (y')} \frac{f(x) }{ |x|^\alpha } dx 
=\frac 1{|y|^\beta } \int_{B_{y_n}^{n-1} (y')} \frac 1{ |x|^\alpha } dx 
\gtrsim y_n^{n-1},
\]
which yields
\[
\int_0^1 \Big( \frac 1{y_n^{\lambda - 1/q}} \frac 1{|y|^\beta } \int_{B_{y_n}^{n-1} (y')} \frac{f(x) }{ |x|^\alpha } dx \Big)^q \frac{dy_n}{y_n }
\gtrsim 
\int_0^1 \frac{dy_n}{y_n^{(\lambda -1/q + 1- n)q +1} }.
\]
However, the integral on the right hand side of the preceding inequality diverges if $$\lambda - 1/q \geq n-1.$$ Hence, we necessarily have $\lambda - 1/q < n-1$. This completes the proof.

\begin{remark}
Now we have some remarks.

\begin{itemize}
 
 \item The identity \eqref{eq-SW-BC-u-n1} reveals that we can transform the HLS inequality on $\R_+^n$ to a suitable HLS inequality on $\R^{n-1}$. This point is fully exploited in \cite{NNN}.
 
 \item Suppose that we only assume $\lambda \geq n-1$, instead of $\lambda \geq n-1/r$. In this scenario, the key estimate \eqref{eq-SW-u-qq} is no longer true because we do not necessarily have $$\alpha + \beta \geq 1/q.$$ 
 
 \item As mentioned in the proof, the above argument works for all $$\lambda \geq n-1/r$$ regardless of $\alpha$ and $\beta$, and this could confuse us why this condition does not depend on $\alpha$ and $\beta$. Nevertheless, seeing \eqref{eq-SW-BC-u-n1} and \eqref{eq-SW-Key-u}, by a simple contradiction argument, we are led to the following inequality for the validity of \eqref{eq-SW-u}
$$\lambda + \alpha + \beta \leq n-1/r,$$
which now depends on $\alpha$ and $\beta$. The above condition tells us that the closer to $n-1/r$ the parameter $\lambda$ is, the smaller the sum $\alpha+\beta$ is.
\end{itemize}
\end{remark}


\subsection{The necessity of (\ref{eq-SW-AB-u}) and (\ref{eq-SW-BC-u})}


In this part we prove the necessity of \eqref{eq-SW-AB-u} and \eqref{eq-SW-BC-u}, whose proof makes use of the same idea used in sections \ref{subsubsec-alpha<}, \ref{subsubsec-beta<}, and \ref{subsubsec-=}. However, to be able to consider the case of $\R^{n-k}\times \R^n$ in Section \ref{sec-Proof-m} below, we provide a quick argument for completeness.

\subsubsection{The necessity of $\beta < n(r-1)/r$}
\label{subsubsec-b-u}

Following the idea used in section \ref{subsubsec-beta<}, to see why $\beta < n(r-1)/r$ is necessary, we simply choose $f = \chi_{B_3^{n-1} \setminus B_2^{n-1}}$ and estimate
\begin{align*}
\int_{\R_+^n} \Big( \int_{\R^{n-1}} & \frac{\chi_{B_3^{n-1} \setminus B_2^{n-1} }(x) dx }{ |x|^\alpha |x-y|^\lambda |y|^\beta } \Big)^q dy
\geq \int_{ B_1^n } 
\Big[ \int_{B_3^{n-1} \setminus B_2^{n-1}} \frac {dx }{ |x|^\alpha |x-y|^\lambda} \Big]^q 
\frac{dy}{|y|^{\beta q}} 
= +\infty
\end{align*}
if $\beta \geq n(r-1)/r=n/q$. Here $q=(1-1/r)^{-1}$ and we need the estimate $$|x-y| \leq \sqrt{2|x|^2 + 3|y|^2} \leq \sqrt{19}$$ to bound the integral of $|x|^{-\alpha}|x-y|^{-\lambda}$ from below away from zero. 


\subsubsection{The necessity of $\alpha < (n-1)(p-1)/p$}
\label{subsubsec-a-u}

Now to see why $\alpha < (n-1)(p-1)/p$ is necessary, we choose $g = \chi_{ B_3^n \setminus B_2^n }$ and estimate
\begin{align*}
\int_{\R^{n-1}} \Big( \int_{\R_+^n} & \frac{\chi_{ B_3^n \setminus B_2^n }(y) dy }{ |x|^\alpha |x-y|^\lambda |y|^\beta } \Big)^q dx
\geq \int_{B_1^{n-1} } 
\Big[ \int_{ B_3^n \setminus B_2^n } \frac {dy }{ |x-y|^\lambda |y|^\beta} \Big]^q 
\frac{dx}{|x|^{\alpha q}} 
= +\infty,
\end{align*}
thanks to $\alpha \geq (n-1)(p-1)/p=(n-1)/q$. Here $q=(1-1/p)^{-1}$ and, as before, we need $$|x-y| \leq \sqrt{2|x|^2 + 3|y|^2} \leq \sqrt{29}$$ to bound the integral of $|x-y|^{-\lambda} |y|^{-\beta}$ from below away from zero. 


\subsubsection{The necessity of $(n-1)/(np) + 1/r + (\lambda + \alpha + \beta + 1)/n =2$}

Still using \eqref{eq-FeGe}, we now obtain
\begin{align*}
\iint_{\R_+^n \times \R^{n-1} } & \frac{f_\epsilon (x) g_\epsilon (y) dx dy}{|x|^\alpha |x-y|^\lambda |y|^\beta } \\
&=\iint_{\R_+^n \times \R^{n-1} } \frac{f (\frac x \epsilon ) g (\frac y \epsilon )dx dy }{\epsilon^{\frac {n-1}p + \frac nr + \lambda + \alpha + \beta} |\frac x\epsilon|^\alpha |\frac x\epsilon - \frac y \epsilon |^\lambda |\frac{y}\epsilon|^\beta } \\
&=\big( \frac 1 \epsilon \big)^{(\frac{n-1}n \frac 1p + \frac 1r + \frac{\lambda + \alpha + \beta + 1}n - 2 ) n}\iint_{\R_+^n \times \R^{n-1} } \frac{f (x ) g (y ) dx dy}{|x|^\alpha |x -y |^\lambda |y |^\beta } .
\end{align*}
From this we obtain the necessity of \eqref{eq-SW-BC-u}.


\subsection{The necessity of (\ref{eq-SW-A+B-u})}
\label{subsec-A+B-u}

As before, to verify the necessity of $\alpha + \beta \geq 0$ for the SW inequality on $\R_+^n$, we construct an example of non-negative functions $f$ and $g$ in such a way that $\|f\|_{L^p(\R^{n-1})} <+\infty$, $\|g\|_{L^r(\R_+^n)} <+\infty$, but
\[
\iint_{\R_+^n \times \R^{n-1}} \frac{f(x)g(y)}{|x|^\alpha |x-y|^\lambda |y|^\beta} dx dy = +\infty.
\] 
To this purpose, we modify our example constructed in subsection \ref{subsec-A+B}. The idea consists of the following two points
\begin{itemize}
 \item to isolate $|x-y|$, as in the previous case, and 
 \item to make sure that $|x|$ and $|y|$ behave similarly, which can be done if we isolate $y_n$.
\end{itemize}
Now fix any small $\epsilon > 0$ and define
\[
f(x)=
\left\{
\begin{aligned}
& (\log_2 x_{n-1})^{-\frac 1p - \epsilon} &&\text{if } \; 2^m \leq x_{n-1} \leq 2^m+1, \; m \geq 1, \\
& & & \text{and } \; |x'| \leq 1,\\
&0 &&\text{otherwise},
\end{aligned}
\right.
\]
and
\[
g(y)=
\left\{
\begin{aligned}
& (\log_2 y_{n-1})^{-\frac 1r - \epsilon} y_n &&\text{if } \; 2^m \leq y_{n-1} \leq 2^m+1, \; m \geq 1,\\
& & & \text{and } \; |y'| \leq 1, \; 1 < y_n < 2,\\
&0 &&\text{otherwise}.
\end{aligned}
\right.
\]
Intuitively, the support of $f$ and $g$ still consists of \textit{disjoint cylinders} concentrated along the half axis $x_{n-1}$. However, to isolate the extra dimension $y_n$, the support of $g$ is chosen as a `slice' of that of $f$. From the above, we obviously have $f \in L^p(\R^{n-1})$. Observe that
\begin{align*}
\int_{\R_+^n} g(y)^r dy &=\sum_{m \geq 1} \int_1^2 \Big( \int_{2^m}^{2^m+1}\Big( \int_{|y''| \leq 1} \frac 1{(\log_2 y_{n-1})^{1+ r\epsilon}} dy' \Big) dy_{n-1} \Big) dy_n\\
&\leq |B_1^{n-2} | \sum_{m \geq 1} \frac{1}{m^{1+ r\epsilon}}
<+\infty,
\end{align*}
which implies $g \in L^r(\R_+^n)$. Now we estimate the double integral as follows
\begin{equation}\label{eq-Double-1-u}
\begin{aligned}
\iint_{\R+^n \times \R^{n-1}} & \frac {f(x)g(y)}{|x|^\alpha |x-y|^\lambda |y|^\beta} dxdy \\
\geq \sum_{m \geq 1}& \iint_{B_1^{n-2} \times [2^m, 2^m+1]} \Big( \iiint_{B_1^{n-2} \times [2^m, 2^m+1] \times [1,2]} \frac {g(y) dy}{ |x-y|^\lambda |y|^\beta} \Big) \frac{f(x) dx}{|x|^\alpha}.
\end{aligned}
\end{equation}
For each $m \geq 1$, let $x$ and $y$ be such that
\begin{equation}\label{eq-Strip-u}
x,(y',y_{n-1}) \in \overline{B_1^{n-2}} \times [2^m, 2^m+1], \quad 1<y_n<2.
\end{equation}
In view of \eqref{eq-Double-1-u}, we need to estimate $|x|^{-\alpha}$, $|y|^{-\beta}$, and $|x-y|^{-\lambda}$ from below. Apparently, we can easily bound $|x-y|$ from the above as follows
\begin{align*}
|x-y|^2 &=|x''-y''|^2 +|x_{n-1}-y_{n-1}|^2 + y_n^2 
\leq 9.
\end{align*}
This together with $\lambda > 0$ allows us to estimate $|x-y|^{-\lambda} \geq 3^{-\lambda}$ provided \eqref{eq-Strip-u} holds for each $m \geq 1$. Since $|x'| \leq 1$ and $x_{n-1} \geq 2$, we also have 
\[
x_{n-1} \leq \sqrt{|x'|^2+x_{n-1}^2} \leq 2x_{n-1}.
\]
Similarly, as $|y'| \leq 1$, $y_{n-1} \geq 2$, and $1 \leq y_n \leq 2$ we get
\[
y_{n-1} \leq \sqrt{|y'|^2+y_{n-1}^2 + y_n^2} \leq 2y_{n-1}.
\]
Hence, regardless of the sign of $\alpha$ and $\beta$, we always have
\[
|x|^{-\alpha} \gtrsim x_{n-1}^{-\alpha}, \quad |y|^{-\beta} \gtrsim y_{n-1}^{-\beta}.
\] 
Putting the above estimates together, we can further estimate \eqref{eq-Double-1-u} as follows
\begin{equation}\label{eq-Double-2-u}
\begin{aligned}
 \iint_{\R_+^n \times \R^{n-1}}& \frac {f(x)g(y)}{|x|^\alpha |x-y|^\lambda |y|^\beta} dxdy\\
\gtrsim& \sum_{m \geq 1} \Big( \iint_{B_1^{n-2} \times [2^m, 2^m+1]} \frac{f(x) dx}{|x|^\alpha} \Big) \Big( \iiint_{B_1^{n-2} \times [2^m, 2^m+1] \times [1,2]} \frac {g(y) dy}{ |y|^\beta} \Big) \\
 =& |B_1^{n-2}|^2 
 \sum_{m \geq 1} \Big( \int_{2^m}^{2^m+1} \frac {dx_{n-1}}{x_{n-1}^\alpha (\log_2 x_{n-1})^{\frac 1p +\epsilon}} \Big) \\
 & \qquad \qquad \times \Big( \int_1^2 y_n dy_n \Big)\Big( \int_{2^m}^{2^m+1} \frac {dy_{n-1}}{y_{n-1}^\beta (\log_2 y_{n-1})^{\frac 1r +\epsilon}} \Big) \\
 \gtrsim& \sum_{m \geq 1} \frac 1{(m+1)^{\frac 1p +\frac 1r +2\epsilon}} 
 \Big( \int_{2^m}^{2^m+1} \frac {dx_{n-1}}{x_{n-1}^\alpha } \Big) \Big( \int_{2^m}^{2^m+1} \frac {dy_{n-1}}{y_{n-1}^\beta } \Big).
\end{aligned}
\end{equation}
Arguing as in section \ref{subsec-A+B}, we know that
\[
\Big( \int_{2^m}^{2^m+1} \frac {dx_{n-1}}{x_{n-1}^\alpha } \Big) \Big( \int_{2^m}^{2^m+1} \frac {dy_{n-1}}{y_{n-1}^\beta } \Big)
\geq \frac 1{(2^m)^{\alpha + \beta}}
\]
provided \eqref{eq-Strip-u} holds. Thus, recalling \eqref{eq-Double-2-u} we arrive at
\begin{align*}
\iint_{\R_+^n \times \R^{n-1}} & \frac {f(x)g(y)}{|x|^\alpha |x-y|^\lambda |y|^\beta} dxdy
\gtrsim \sum_{m \geq 1} \frac 1{(2^m)^{\alpha + \beta}} \frac 1{(m+1)^{\frac 1p +\frac 1r + 2\epsilon}} .
\end{align*}
From this we necessarily have $\alpha +\beta \geq 0$. This proves the necessity of \eqref{eq-SW-AB-u}.


\subsection{The necessity of (\ref{eq-SW-Key-u})}
\label{subsec-P+R-u}

Now we establish the necessity of the condition $1/p+1/r \geq 1$. The example provided below is essentially the same as that constructed in section \ref{subsec-P+R} with some necessary changes due to the fact that we now work on $\R_+^n \times \R^{n-1}$. To this purpose, by duality, we obtain an equivalent form of \eqref{eq-SW-u} as follows
\begin{equation}\label{eq-SW-u-n2}
\int_{\R^{n-1}} \Big( \int_{\R_+^n} \frac{g (y) }{|x|^\alpha |x-y|^\lambda |y|^{\beta }} dy \Big)^q dx
\lesssim \| g\|_{L^p(\R_+^n)}
\end{equation}
with
\begin{equation}\label{eq-SW-u-qq}
q = \Big(1- \frac 1p \Big)^{-1} \overset{\eqref{eq-SW-BC-u-n1}}{=}
\Big(\frac 1 r + \frac {\lambda + \alpha + \beta - \frac 1q }{n-1} -1 \Big)^{-1}
<r,
\end{equation}
thanks to $\lambda + \alpha + \beta - 1/q >n-1$. Now we choose 
\[
g(y)=
\left\{
\begin{aligned}
& |y|^{-\frac nr } (\log |y|)^{- \frac 1 q} & & \text{if } |y| \geq 2,\\
& 0 & & \text{otherwise}.
\end{aligned}
\right.
\]
As in subsection \ref{subsec-P+R}, we clearly have
\[
\int_{\R^n} g(y )^r dy = |\mathbb S^{n-1}| \int_2^{+\infty} \frac{d\rho}{\rho (\log \rho)^{ \frac r q}} < +\infty,
\]
thanks to $r > q$. Hence, $ g \in L^r (\R_+^n)$. Clearly, by \eqref{eq-SW-u-n2} we have
\begin{align*}
+\infty > \| g\|_{L^p(\R_+^n)}^q &> \int_{\R^{n-1}} \Big( \int_{\R_+^n} \frac{g (y) }{|x|^\alpha |x-y|^\lambda |y|^\beta} dy \Big)^q dx\\
&\geq
\int_{|x| \geq 4} \Big( \int_{|x|/2 \leq |y| \leq 2|x|} \frac{ dy }{|y|^{\beta + \frac nr} (\log |y|)^{\frac 1 q } |x-y|^\lambda } \Big)^q \frac{ dx}{|x|^{\alpha q} }.
\end{align*}
As $\lambda > 0$ and $|y| \leq 2|x|$, we can bound $|x-y|$ from the above as follows
\[
|x-y|^2 = |x-y'|^2 + y_n^2 \leq 2|x|^2 + 2|y'|^2 + y_n^2 
\leq 6|x|^2,
\]
which, together with $\lambda>0$, implies
\[
|y|^{\beta + \frac nr} |x-y|^{\lambda} \lesssim |x|^{\lambda +\beta + \frac n r}
\]
for any $|x|/2 \leq |y| \leq 2|x|$. Hence, together with 
\[
0< \log |y| \leq \log (2|x|) \leq \log 2 + \log |x| \leq 2\log |x|,
\] 
we know that
\[
\frac 1{|y|^{\beta + \frac nr} (\log |y|)^{\frac 1 q } |x-y|^\lambda }
\gtrsim \frac 1{|x|^{ \lambda + \beta + \frac nr } (\log |x|)^{ \frac 1 q } }.
\]
Observe that
\begin{equation}\label{eq-SW-u-Exponent}
\begin{aligned}
\Big(\lambda + \alpha +\beta + \frac nr \Big) q 
= \Big(\frac{n-1}n \big( 1 - \frac 1p \big) + 1 \Big) n q = n-1 + n q.
\end{aligned}
\end{equation}
Hence
\begin{align*}
+\infty > \| g \|_{L^p(\R_+^n )}^q 
&\gtrsim
\int_{|x| \geq 4} \Big( \int_{|x|/2 \leq |y| \leq 2|x|} dy \Big)^q
\frac {dx}{|x|^{( \lambda + \alpha +\beta + \frac nr )q } (\log |x|) } \\
&=
\int_{|x| \geq 4} \frac 1{|x|^{n-1} \log |x| } 
= +\infty,
\end{align*}
thanks to $\int_{|x|/2 \leq |y| \leq 2|x|} dy \sim |x|^n$. Hence, we obtain the necessity of \eqref{eq-SW-Key-u} as claimed.


\section{Necessity of conditions for (\ref{eq-SW-m}) on $\R^{n-k} \times \R^n$: proof of Theorem \ref{thm-MAIN-m}}
\label{sec-Proof-m}

In this section, we show that our argument performed in section \ref{sec-Proof-u} above can be further applied to the case of $\R^{n-k} \times \R^n$ with some changes due to the fact that we are no longer in $\R^{n-1} \times (0, +\infty)$. The key idea is to transform the higher-dimensional space $\R^k$ into the one-dimensional case $(0,+\infty)$. We also note that if we let $k=0$, then Theorem \ref{thm-MAIN-m} becomes Theorem \ref{thm-MAIN}. However, to further illustrate the method, we shall use the condition $k \geq 1$ in several places.

For clarity and convenience, for a general point $x \in \R^{n-k}$, we write 
\[
x=(x',x_{n-k}) \in \R^{n-k-1} \times \R.
\] 
Hence, for $y \in \R^n$ we shall write 
\[
y = (y', y_{n-k}, y'') \in \R^{n-k-1} \times \R \times \R^k.
\] 


\subsection{The necessity of (\ref{eq-SW-L-m})}
\label{subsubsec-L-m}

Employing a similar argument as in subsection \ref{subsec-L-u}, we can easily rule out the case $\lambda \leq 0$. Now we rule out the case $\lambda \geq n-k/r$. For some non-negative function $f \in L^p(\R^{n-k})$ to be determined later, the inequality \eqref{eq-SW-BC-m} is equivalent to
\begin{align*}
\| f \|_{L^p(\R^{n-k})}^q &\gtrsim \int_{\R^n} \Big( \int_{\R^{n-k}} \frac{f(x) }{ |x|^\alpha |x-y|^\lambda |y|^\beta } dx \Big)^q dy \\
&=|\mathbb S^{k-1}| \int_{\R^{n-k}}\Big[ \int_0^{+\infty} \Big( \frac 1{ |y|^\beta } \int_{\R^{n-k}} \frac{f(x) }{ |x|^\alpha |x-y|^\lambda } dx \Big)^q \rho^{k-1} d\rho \Big] dy',
\end{align*}
with $q=r/(r-1)$ and $\rho = |y''|$. Now we need some estimate for the integral $\int_{\R^{n-k}}$ which costs us some energy due to the fact that $k \geq 1$. Still by Fubini's theorem and the non-decreasing of 
$$ \rho \mapsto \int_{B_{\rho}^{n-k} (y',y_{n-k})} |x|^{-\alpha } f(x)dx,$$
we easily get
\begin{align*}
\int_{\R^{n-k}} \frac{f(x) }{ |x|^\alpha |x-y|^\lambda } dx
&=\lambda \int_{\R^{n-k}} \Big( \int_{|x-y|}^{+\infty} \frac{d\rho}{\rho^{\lambda +1}} \Big) \frac{f(x) }{ |x|^\alpha } dx\\
&=\lambda \int_{|y''|}^{+\infty} \Big( \int_{B_{\rho - |y''|}^{n-k} (y',y_{n-k})} \frac{f(x) }{ |x|^\alpha } dx \Big) \frac{d\rho}{\rho^{\lambda +1}}\\
&\gtrsim \int_{2|y''|}^{+\infty} \Big( \int_{B_{\rho - |y''|}^{n-k} (y',y_{n-k})} \frac{f(x) }{ |x|^\alpha } dx \Big) \frac{d\rho}{\rho^{\lambda +1}}\\
& \geq \Big( \int_{B_{|y''|}^{n-k} (y',y_{n-k})} \frac{f(x) }{ |x|^\alpha } dx \Big) \int_{2|y''|}^{+\infty} \frac{d\rho}{\rho^{\lambda +1}}\\
&\gtrsim \frac 1{|y''|^\lambda } \int_{B_{|y''|}^{n-k} (y',y_{n-k})} \frac{f(x) }{ |x|^\alpha } dx.
\end{align*}
Hence, we obtain
\begin{align*}
\| f \|_{L^p(\R^{n-k})}^q 
&\gtrsim |\mathbb S^{k-1}| \int_{\R^{n-k}} \Big[\int_0^{+\infty} \Big( \frac 1{ |y|^{\beta } } \int_{\R^{n-k}} \frac{f(x) }{ |x|^\alpha |x-y|^\lambda} dx \Big)^q \rho^{k-1} d\rho \Big] dy' dy_{n-k}\\
&\gtrsim \int_{B_4^{n-k} \setminus B_2^{n-k}} \Big[ \int_0^1 \Big( \frac 1{\rho^{ \lambda - k/q} |y|^\beta } \int_{B_{\rho}^{n-k} (y', y_{n-k})} \frac{f(x) }{ |x|^\alpha } dx \Big)^q \frac{d\rho}{\rho} \Big] dy' dy_{n-k} .
\end{align*}
For the last line in the above computation, thanks to $|(y', y_{n-k})|^2 \geq 4$ and $0\leq |y''| \leq 1$, we know that
\[
B_1^{n-k} \subset B_{\rho}^{n-k} (y', y_{n-k}) \subset B_5^{n-k}
\]
and that $2 \leq |y| \leq \sqrt 5$. Hence, if we choose $f = \chi_{B_6^{n-1}}$, then we can bound
\[
\frac 1{|y|^\beta } \int_{B_{\rho}^{n-k} (y', y_{n-k})} \frac{f(x) }{ |x|^\alpha } dx 
 \gtrsim \rho^{n-k},
\]
which yields
\[
\int_0^1 \Big( \frac 1{\rho^{\lambda - k/q}} \frac 1{|y|^\beta } \int_{B_{\rho}^{n-k} (y', y_{n-k})} \frac{f(x) }{ |x|^\alpha } dx \Big)^q \frac{d \rho}{\rho }
\gtrsim 
\int_0^1 \frac{d\rho}{\rho^{(\lambda -k/q -(n- k))q +1} }.
\]
However, the integral on the right hand side of the preceding inequality diverges if
\[
\lambda - k/q \geq n-k.
\] 
Hence, we necessarily have $\lambda - k/q < n-k$. This completes the proof.

\begin{remark}
The necessity of the condition $\lambda - k/q < n-k$ for \eqref{eq-SW-m} is quite different from that of \eqref{eq-SW-m-simpler} as \eqref{eq-SW-m-simpler} holds for any $\lambda > 0$; see \cite{NNN}.
\end{remark}


\subsection{The necessity of (\ref{eq-SW-BC-m}) and (\ref{eq-SW-AB-m})}

As before, the balance condition \eqref{eq-SW-BC-m} follows from standard computation. In addition, the condition \eqref{eq-SW-AB-m} can also be proved by using the idea in subsections \ref{subsubsec-a-u} and \ref{subsubsec-b-u}. Hence we omit the details.


\subsection{The necessity of (\ref{eq-SW-A+B-m})}
\label{subsec-A+B-m}

To verify the necessity of $\alpha + \beta \geq 0$ for the inequality \eqref{eq-SW-m} on $\R^{n-k} \times \R^n$, we follows the construction in subsection \ref{subsec-A+B-u}. Now fix any small $\epsilon > 0$ and define
\[
f(x)=
\left\{
\begin{aligned}
& (\log_2 x_{n-k})^{-\frac 1p - \epsilon} &&\text{if } \; 2^m \leq x_{n-k} \leq 2^m+1, \; m \geq 1, \\
& & & \text{and } \; |x'| \leq 1,\\
&0 &&\text{otherwise},
\end{aligned}
\right.
\]
and
\[
g(y)=
\left\{
\begin{aligned}
& (\log_2 y_{n-k})^{-\frac 1r - \epsilon} |y''| &&\text{if } \; 2^m \leq y_{n-k} \leq 2^m+1, \; m \geq 1,\\
& & & \text{and } \; |y'| \leq 1, \; 1 < |y''| < 2,\\
&0 &&\text{otherwise}.
\end{aligned}
\right.
\]
From the above construction, it is clear that the quantity $|y''|$ serves a similar role as $y_n$ in the case of $\R^{n-1} \times \R_+^n$. Now we obviously have $f \in L^p(\R^{n-k})$. Recall that $\int_{\R^k} = |\mathbb S^{k-1}| \int_0^{+\infty}$. Then for the function $g$ we have
\begin{align*}
\int_{\R^n} g(y)^r dy &=|\mathbb S^{k-1}| \sum_{m \geq 1} \int_1^2 \Big( \int_{2^m}^{2^m+1}\Big( \int_{|y'| \leq 1} \frac 1{(\log_2 y_{n-k})^{1+ r\epsilon}} dy' \Big) dy_{n-k} \Big) d\rho\\
&\leq |\mathbb S^{k-1}| |B_1^{n-k-1} | \sum_{m \geq 1} \frac{1}{m^{1+ r\epsilon}}
<+\infty,
\end{align*}
which implies $g \in L^r(\R^n)$. Now we estimate the double integral as follows
\begin{equation}\label{eq-Double-1-m}
\begin{aligned}
\iint_{\R^n \times \R^{n-k}} & \frac {f(x)g(y)}{|x|^\alpha |x-y|^\lambda |y|^\beta} dxdy \\
\geq \sum_{m \geq 1}& \iint_{B_1^{n-k-1} \times [2^m, 2^m+1] } \Big( \iiint_{B_1^{n-k-1} \times [2^m, 2^m+1] \times [1, 2]} \frac {g(y) dy}{ |x-y|^\lambda |y|^\beta} \Big) \frac{f(x) dx}{|x|^\alpha}.
\end{aligned}
\end{equation}
For each $m \geq 1$, let $x$ and $y$ be such that
\begin{equation*}\label{eq-Strip-m}
x, (y',y_{n-k}) \in \overline{B_1^{n-k-1}} \times [2^m, 2^m+1], \quad 1<|y''|<2.
\end{equation*}
As before, we need to estimate $|x|^{-\alpha}$, $|y|^{-\beta}$, and $|x-y|^{-\lambda}$ from below. Apparently, we can easily bound $|x-y|$ from the above as follows
\begin{align*}
|x-y|^2 &=|x'-y'|^2 +|x_{n-k}-y_{n-k}|^2 + |y''|^2 
\leq 9,
\end{align*}
giving $|x-y|^{-\lambda} \geq 3^{-\lambda}$. Now as $|x'| \leq 1$ and $x_{n-k} \geq 2$, we also have $x_{n-k} \leq |x \leq 2x_{n-k}$, giving $$|x|^{-\alpha} \gtrsim x_{n-k}^{-\alpha}.$$ Similarly, as $|y'| \leq 1$, $y_{n-k} \geq 2$, and $1 \leq |y''| \leq 2$ we get $y_{n-k} \leq |y|\leq 2y_{n-k}$, giving $$|y|^{-\beta} \gtrsim y_{n-k}^{-\beta}.$$ Putting the above estimates together and similar to \eqref{eq-Double-2-u} we can further estimate \eqref{eq-Double-1-m} as follows
\begin{equation*}\label{eq-Double-2-m}
\begin{aligned}
 \iint_{\R^n \times \R^{n-k}}& \frac {f(x)g(y)}{|x|^\alpha |x-y|^\lambda |y|^\beta} dxdy\\
\gtrsim& \sum_{m \geq 1} \Big( \int_{2^m}^{2^m+1} \frac {dx_{n-k}}{x_{n-k}^\alpha (\log_2 x_{n-k})^{\frac 1p +\epsilon}} \Big) \\
 & \qquad \times \Big( |\mathbb S^{k-1}| \int_1^2 \rho^k d \rho \Big)\Big( \int_{2^m}^{2^m+1} \frac {dy_{n-k}}{y_{n-k}^\beta (\log_2 y_{n-k})^{\frac 1r +\epsilon}} \Big) \\
 \gtrsim& \sum_{m \geq 1} \frac 1{(m+1)^{\frac 1p +\frac 1r +2\epsilon}} 
 \Big( \int_{2^m}^{2^m+1} \frac {dx_{n-k}}{x_{n-k}^\alpha } \Big) \Big( \int_{2^m}^{2^m+1} \frac {dy_{n-k}}{y_{n-k}^\beta } \Big) . 
\end{aligned}
\end{equation*}
Thus, we arrive at
\begin{align*}
\iint_{\R^n \times \R^{n-k}} & \frac {f(x)g(y)}{|x|^\alpha |x-y|^\lambda |y|^\beta} dxdy
\gtrsim \sum_{m \geq 1} \frac 1{(2^m)^{\alpha + \beta}} \frac 1{(m+1)^{\frac 1p +\frac 1r + 2\epsilon}} .
\end{align*}
This proves the necessity of \eqref{eq-SW-AB-m}, thanks to \eqref{eq-SW-m}.


\subsection{The necessity of (\ref{eq-SW-Key-m})}
\label{subsec-P+R-m}

Now we establish the necessity of the condition $1/p+1/r \geq 1$ by following the construction in subsection \ref{subsec-P+R-u} with some necessary changes due to the fact that we now work on $\R^{n-k} \times \R^n$. To this purpose, by duality, we obtain an equivalent form of \eqref{eq-SW-m} as follows
\begin{equation}\label{eq-SW-m-n2}
\int_{\R^{n-k}} \Big( \int_{\R^n} \frac{g (y) }{|x|^\alpha |x-y|^\lambda |y|^{\beta }} dy \Big)^q dx
\lesssim \| g\|_{L^p(\R^n)}^q
\end{equation}
with $q = (1- 1/p)^{-1} < r$, thanks to $\lambda + \alpha + \beta - 1/q >n-k$. Now we choose 
\[
g(y)=
\left\{
\begin{aligned}
& |y|^{-\frac nr } (\log |y|)^{- \frac 1 q} & & \text{if } |y| \geq 2,\\
& 0 & & \text{otherwise}.
\end{aligned}
\right.
\]
We clearly have $ g \in L^r (\R^n)$, thanks to $r > q$. Clearly, by \eqref{eq-SW-u-n2} we have
\begin{align*}
+\infty > \| g\|_{L^p(\R^n)}^q &> \int_{\R^{n-k}} \Big( \int_{\R^n} \frac{g (y) }{|x|^\alpha |x-y|^\lambda |y|^\beta} dy \Big)^q dx\\
&\geq
\int_{|x| \geq 4} \Big( \int_{|x|/2 \leq |y| \leq 2|x|} \frac{ dy }{|y|^{\beta + \frac nr} (\log |y|)^{\frac 1 q } |x-y|^\lambda } \Big)^q \frac{ dx}{|x|^{\alpha q} }.
\end{align*}
As $\lambda > 0$ and $|y| \leq 2|x|$, we can bound $|x-y|$ from the above as follows
\begin{align*}
|x-y|^2 &= |x'-y'|^2 + |x_{n-k} - x_{n-k}|^2 + |y''|^2 \\
& \leq 2|x|^2 + 2|y'|^2 + 2 y_n^2 + |y''|^2 \\
&\leq 6|x|^2,
\end{align*}
which implies
\[
|y|^{\beta + \frac nr} |x-y|^{\lambda} \lesssim |x|^{\lambda +\beta + \frac n r}
\]
for any $|x|/2 \leq |y| \leq 2|x|$. Hence, we know that
\[
\frac 1{|y|^{\beta + \frac nr} (\log |y|)^{\frac 1 q } |x-y|^\lambda }
\gtrsim \frac 1{|x|^{ \lambda + \beta + \frac nr } (\log |x|)^{ \frac 1 q } }.
\]
Hence, together with \eqref{eq-SW-m-n2} we can estimate
\begin{align*}
+\infty > \| g \|_{L^p(\R^n )}^q
&\gtrsim
\int_{|x| \geq 4} \Big( \int_{|x|/2 \leq |y| \leq 2|x|} dy \Big)^q
\frac {dx}{|x|^{( \lambda + \alpha +\beta + \frac nr )q } (\log |x|) } \\
&=
\int_{|x| \geq 4} \frac 1{|x|^{n-k} \log |x| } = +\infty,
\end{align*}
thanks to $ (\lambda + \alpha +\beta + n/r ) q = n-k + n q$. This proves the necessity of \eqref{eq-SW-Key-m} as claimed.


\section*{Acknowledgments}

This project was initiated while the author was staying at the University of Tokyo. Funding received including the Tosio Kato Fellowship awarded in 2018 is greatly acknowledged. The author also benefited from the Vietnam Institute for Advanced Study in Mathematics (VIASM) during his visit in 2021 and from the Vietnam National University, Hanoi (VNU) under project number QG.21.01. Last but not least, we thank Jingbo Dou and Quoc Hung Nguyen for useful discussion during the preparation of the note.




\begin{thebibliography}{9999999}

\bibitem[CLLT18]{CLLT-TAMS}
\textsc{L. Chen, Z. Liu, G. Lu, and C. Tao},
Reverse Stein--Weiss inequalities and existence of their extremal functions,
\textit{Trans. Amer. Math. Soc.} \textbf{370} (2018) 8429--8450.

\bibitem[CLLT20]{CLLT}
\textsc{L. Chen, Z. Liu, G. Lu, and C. Tao},
Stein--Weiss inequalities with the fractional Poisson kernel,
\textit{Rev. Mat. Iberoam.} \textbf{36} (2020) 1289--1308.

\bibitem[CLT19]{CLT-ANS}
\textsc{L. Chen, G. Lu, and C. Tao},
Reverse Stein--Weiss inequalities on the upper half space and the existence of their extremals,
\textit{Adv. Nonlinear Stud.} \textbf{19} (2019) 475--494.

\bibitem[Dou16]{Dou2016}
\textsc{J. Dou},
Weighted Hardy--Littlewood--Sobolev inequalities on the upper half space,
\textit{Commun. Contemp. Math.} \textbf{18} (2016) Art. 1550067.

\bibitem[DZ15]{DZ-IMRN}
\textsc{J. Dou and M. Zhu},
Sharp Hardy--Littlewood--Sobolev inequality on the upper half space,
\textit{Int. Math. Res. Not. IMRN} \textbf{2015} (2015) 651--687.

\bibitem[HLZ12]{HLZ12}
\textsc{X. Han, G. Lu, and Z. Zhu},
Hardy-Littlewood-Sobolev and Stein-Weiss inequalities and integral systems on the Heisenberg group,
\textit{Nonlinear Anal.} \textbf{75} (2012) 4296--4314. 

\bibitem[HWY08]{HWY08}
\textsc{F. Hang, X.D. Wang, and X.D. Yan},
Sharp integral inequalities for harmonic functions, 
\textit{Comm. Pure Appl. Math.} \textbf{61} (2008) 54--95.

\bibitem[HL28]{hl1928}
\textsc{G.H. Hardy and J.E. Littlewood},
Some properties of fractional integrals. I,
\textit{Math. Z.} \textbf{27} (1928) 565-606.


\bibitem[KRS19]{KRS}
\textsc{A. Kassymov, M. Ruzhansky, and D. Suragan},
Hardy--Littlewood--Sobolev and Stein--Weiss inequalities on homogeneous Lie groups,
\textit{Integral Transforms Spec. Funct.} \textbf{30} (2019) 643--655. 

\bibitem[Lie83]{Lieb}
\textsc{E.H. Lieb},
Sharp constants in the Hardy--Littlewood--Sobolev and related inequalities,
\textit{Ann. of Math.} \textbf{118} (1983) 349--374.

\bibitem[Lio84]{Lions}
\textsc{P.L. Lions},
The concentration-compactness principle in the calculus of variations. The locally
compact case. I, 
\textit{Ann. Inst.H. Poincar\'e Anal. Non Lin\'eaire} \textbf{1} (1984) 109--145.

\bibitem[NNN20]{NNN}
\textsc{Q.A. Ng\^o, Q.H. Nguyen, and V.H. Nguyen},
An optimal Hardy--Littlewood--Sobolev inequality on $\R^{n-k} \times \R^n$ and its consequences, arXiv:2009.09868.





\bibitem[SDY14]{SDY2014}
\textsc{Z. Shi, W. Di, and D. Yan},
Necessary and sufficient conditions of doubly weighted Hardy--Littlewood--Sobolev inequality,
\textit{Anal. Theory Appl.} \textbf{30} (2014) 193--204.

\bibitem[Sob38]{sobolev1938}
\textsc{S.L. Sobolev},
On a theorem of functional analysis,
\textit{Math. Sb. (N.S.)} {\bf 4} (1938) 471-479. English transl. in \textit{Amer. Math. Soc. Transl. Ser. 2} {\bf 34} (1963) 39-68.

\bibitem[SW58]{SW}
\textsc{E.M. Stein and G. Weiss},
Fractional integrals in $n$-dimensional Euclidean space,
\textit{J. Math. Mech.} {\bf 7} (1958) 503-514.

\end{thebibliography}
\end{document}